\documentclass[12pt]{article}

\usepackage{amssymb,amsmath,latexsym}
\usepackage[active]{srcltx}
\usepackage[matrix,arrow,curve]{xy}
\input xy
\xyoption{all}

\oddsidemargin=\evensidemargin \advance\topmargin by-\headheight

\DeclareMathOperator{\aut}{Aut}

\DeclareMathOperator{\cyc}{Cyc}

\DeclareMathOperator{\Ext}{Ext}

\DeclareMathOperator{\id}{id}
\DeclareMathOperator{\hol}{Hol}

\DeclareMathOperator{\rad}{rad}
\DeclareMathOperator{\rk}{rk}

\DeclareMathOperator{\sym}{Sym}

\def\Q{{\mathbb Q}}

\def\mZ{{\mathbb Z}}

\def\A{{\cal A}}

\def\H{{\cal G}}

\def\M{{\cal M}}

\def\S{{\cal S}}

\def\twoe{\underset{\scriptscriptstyle ^2}{\approx}}

\def\lg{\langle}
\def\ov{\overline}
\def\rg{\rangle}

\def\wt{\widetilde}

\def\proof{{\bf Proof}.\ }
\def\bull{\vrule height .9ex width .8ex depth -.1ex }

\def\qaq{\quad\text{and}\quad}
\def\VRT#1{*=<6mm>[o][F-]{#1}}

\makeatletter
\renewcommand{\subsection}{\@startsection{subsection}{2}{0mm}{-2mm}{-2mm}
{\bf\normalsize}}

\makeatother

\newtheorem{formula}{}[section]
\newtheorem{proposition}[formula]{Proposition}
\newtheorem{definition}[formula]{Definition}
\newtheorem{corollary}[formula]{Corollary}
\newtheorem{remark}[formula]{Remark}
\newtheorem{lemma}[formula]{Lemma}
\newtheorem{theorem}[formula]{Theorem}

\newtheorem{example}[formula]{Example}

\def\thrm{\begin{theorem}}
\def\thrml#1{\begin{theorem}\label{#1}}
\def\ethrm{\end{theorem}}
\def\prpstn{\begin{proposition}}
\def\prpstnl#1{\begin{proposition}\label{#1}}
\def\eprpstn{\end{proposition}}
\def\rmrk{\begin{remark}}
\def\rmrkl#1{\begin{remark}\label{#1}}
\def\ermrk{\end{remark}}
\def\dfntn{\begin{definition}}
\def\dfntnl#1{\begin{definition}\label{#1}}
\def\edfntn{\end{definition}}
\def\nmrt{\begin{enumerate}}
\def\enmrt{\end{enumerate}}
\def\tm#1{\item[{\rm (#1)}]}
\def\qtn{\begin{equation}}
\def\qtnl#1{\begin{equation}\label{#1}}
\def\eqtn{\end{equation}}
\def\lmm{\begin{lemma}}
\def\lmml#1{\begin{lemma}\label{#1}}
\def\elmm{\end{lemma}}
\def\crllr{\begin{corollary}}
\def\crllrl#1{\begin{corollary}\label{#1}}
\def\ecrllr{\end{corollary}}
\def\hpthss{\begin{hypothesis}}
\def\hpthssl#1{\begin{hypothesis}\label{#1}}
\def\ehpthss{\end{hypotxesis}}
\def\xmpl{\begin{example}}
\def\xmpll#1{\begin{example}\label{#1}}
\def\exmpl{\end{example}}
\def\css{\begin{cases}}
\def\ecss{\end{cases}}

\begin{document}

\title{Characterization of cyclic Schur groups}
\author{
Sergei Evdokimov \\[-2pt]
\small Steklov Institute of Mathematics\\[-4pt]
\small at St. Petersburg \\[-4pt]
{\tt \small evdokim@pdmi.ras.ru }
%\thanks{The work was partially supported by RFFI Grant 07-01-00485.}
\thanks{The work was partially supported by
Slovenian-Russian bilateral project, grant no. BI-RU/10-11-018.}
\and
Istv\'an Kov\'acs\\[-2pt]
\small University of Primorska, Koper, Slovenia\\[-4pt]
{\tt \small kovacs@pef.upr.si}
\thanks{The work was partially supported by
Slovenian-Russian bilateral project, grant no. BI-RU/10-11-018.}
\and
Ilya Ponomarenko\\[-2pt]
\small Steklov Institute of Mathematics\\[-4pt]
\small at St. Petersburg \\[-4pt]
{\tt \small inp@pdmi.ras.ru}
\thanks{The work was partially supported by
Slovenian-Russian bilateral project, grant no. BI-RU/10-11-018
and RFFI Grant 11-01-00760-a.}
}

\date{}
\maketitle

\begin{abstract}
A finite group $G$ is called a Schur group, if any Schur ring over~$G$ is
associated in a natural way with a subgroup of $\sym(G)$ that contains all right 
translations. It was proved by R.~P\"oschel (1974) that given a prime
$p\ge 5$ a $p$-group is Schur if and only if it is cyclic. We prove that
a cyclic group of order $n$ is Schur if and only if $n$ belongs to one
of the following five families of integers:
$p^k$, $pq^k$, $2pq^k$, $pqr$, $2pqr$ where $p,q,r$ are distinct primes,
and $k\ge 0$ is an integer.
\end{abstract}

\section{Introduction}\label{240111b}
Let $G$ be a finite group. A subring of the group ring $\Q G$ is called a {\it Schur ring}
or {\it S-ring} over~$G$, if it closed with respect to the componentwise multiplication
and inversion. The first construction of such a ring was proposed by I.~Schur~\cite{S33}
in connection with his famous result on permutation groups containing a
regular cyclic subgroup. Namely, let $\Gamma$ be a permutation group on the set~$G$ that contains
the regular group $G_{right}$ induced by right multiplications,
$$
G_{right}\le\Gamma\le \sym(G).
$$
Denote by $\Gamma_1$ the stabilizer of the identity of $G$ in $\Gamma$.
Then the  submodule of~$\Q G$ spanned by the $\Gamma_1$-orbits 
(transitivity module) is an S-ring over~$G$. Such an S-ring was called {\it schurian}
in~\cite{P74}. The general theory of S-rings was developed by
H.~Wielandt in~\cite{W64} where in particular he constructed an S-ring
which cannot be obtained by the Schur method.\medskip
%. In Theorem~26.4 of that book 

\noindent {\bf Definition (R.~P\"oschel).} {\it A finite group $G$ is called Schur,
if any S-ring over~$G$ is schurian.}\medskip

The Wielandt example shows that not every finite group is Schur. More exactly, he
proved that the group $\mZ_p\times\mZ_p$ is not Schur for prime $p\ge 5$. This
fact was used by R.~P\"oschel in~\cite{P74} to prove the following theorem.\medskip

\noindent {\bf Theorem.} {\it Any section of a Schur group is a Schur
group. Moreover, for a prime $p\ge 5$ a $p$-group is Schur if and only if
it is cyclic.\bull}\medskip

Since any finite nilpotent group is a direct product of its Sylow subgroups, we
immediately obtain the following result.\medskip

\noindent {\bf Corollary.} {\it A nilpotent group of order coprime to~$6$ is Schur
only if it is cyclic.\bull}\medskip

The above results show the importance of the cyclic case for
the characterization of Schur groups. It should be noted that
by the P\"oschel theorem any cyclic $p$-group is Schur for $p\ge 5$.
In fact, the schurity of cyclic $3$-groups was also proved in~\cite{P74},
whereas the same result for $p=2$ was obtained in~\cite{K88}. However,
till 2001 no cyclic non-Schur group was known,
and moreover it was conjectured that all cyclic groups are Schur (the Schur-Klin
conjecture). This conjecture had also been supported by the fact that
the group $\mZ_n$ where $n$ is a product of two distinct primes,
is a Schur one~\cite{KP}. The first counterexamples to the conjecture
were constructed in~\cite{EP01ae}; in all these examples $n$ was the product of
at least four primes. Later in~\cite{EP11} the schurity of $\mZ_n$
was proved when $n$ is the product of at most three primes or
$n=p^3q$ where $p$ and $q$ are distinct primes. The main result of this paper
completes the characterization of cyclic Schur groups.

\thrml{040211a}
A cyclic group of order $n$ is Schur if and only if $n$ belongs to one
of the following five (partially overlapped) families of integers:
\qtnl{250311b}
p^k,\ pq^k,\ 2pq^k,\ pqr,\ 2pqr
\eqtn
where $p,q,r$ are distinct primes, and $k\ge 0$ is an integer.
\ethrm

\crllrl{160911a}
The minimum order of a cyclic non-Schur group equals~$72$.\bull
\ecrllr

Let us briefly outline the proof of Theorem~\ref{040211a}.
To prove the necessity, for each integer~$n$ satisfying the hypothesis of
Theorem~\ref{040211b} we construct explicitly a non-schurian S-ring over a group $\mZ_n$.
This ring is the generalized wreath
product of two smaller schurian S-rings each of which is in its turn the generalized
wreath product of normal S-rings; the way is essentially the same as one used
in~\cite{EP01ae}. It turns out (Lemma~\ref{250311a}) that the complement to the set of all  these $n$
coincides with the set of all numbers listed in~\eqref{250311b}.\medskip

%It was mentioned above that any cyclic group the order of which belongs to the first or
%fourth family is a Schur one. Therefore the proof of the sufficiency is reduced
%to proving the schurity of any S-ring over $\mZ_n$ where $n$ belongs to the fifth,
%second or third families. Assuming on the contrary that there
%is a non-schurian S-ring, we successively come to a contradiction
%for each of these cases in Sections~\ref{250711c},
%\ref{250711a} and~\ref{250711b} respectively. 

To prove the sufficiency %of Theorem~\ref{040211a}
we have to verify that any S-ring over a cyclic group of order $n$ belonging to one of 
families~\eqref{250311b}, is schurian. We observe that any divisor of such $n$ also belongs to at least
one of these families.

\dfntn
A non-schurian S-ring $\A$ over a group $G$ is called minimal if the S-ring $\A_S$ is schurian
for any $\A$-section $S\ne G/1$.
\edfntn

It is easily seen that any non-schurian S-ring contains a section the restriction to which is minimal 
non-schurian. Thus the sufficiency in Theorem~\ref{040211a} immediately follows from the theorem below.

\thrml{130112a}
The order of the underlying group of a minimal non-schurian circulant S-ring cannot belong to any of 
families~\eqref{250311b}.
\ethrm  

There are two key observations to prove Theorem~\ref{130112a} that are based on the results of~\cite{EP11}. 
The first is that any non-schurian circulant S-ring %in question can be assumed 
is a fusion of a quasidense~\footnote{Quasidense circulant 
S-rings are introduced and studied in Section~\ref{130712a}.} non-schurian circulant S-ring (Theorem~\ref{220911a}).
The second is that such an S-ring is a proper generalized wreath product of two smaller schurian quasidense 
S-rings. % (Lemma~\ref{160112a}). Moreover, 
Moreover, in the minimal case each of them
is in its turn a proper generalized wreath product (Theorem~\ref{160112b}). %A further analysis in 
We use these observations in Sections~\ref{190312a} and~\ref{160312a}  to exclude the first, second and 
fourth families, and the
case $p=2$ in the other two families. The proof is completed in Section~\ref{190312b} by applying
the criterion of schurity for S-rings of special form that was proved 
in Section~\ref{130712c}. It should be mentioned that throughout the proof of Theorem~\ref{130112a} we 
use several auxiliary results on circulant S-rings that are collected in Section~\ref{130712b}.\medskip

%Sections~\ref{230911a} and~\ref{230911b} contain auxiliary lemmas which
%are used to prove the sufficiency.

In this paper we follow the notation and terminology of paper~\cite{EP11}. When referring to this paper
we keep only the number of the statement, preceding it by the letter~A (e.g. instead of \cite[Theorem~4.1]{EP11}
we write Theorem~A4.1). Several additional notations are listed below.\medskip\medskip

%{\bf Additional notation.}

We write $\A\cong\A'$ when S-rings $\A$ and $\A'$ are Cayley isomorphic.\medskip

For an S-ring $\A$ and an $\A$-section~$S$ we set $\hol_\A(S)=\hol(S)\cap\aut(\A_S)$.\medskip

For an S-ring $\A$ over a group $G$ we set 
$$
\M(\A)=\{\Gamma\le\aut(\A):\ \Gamma\twoe\aut(\A)\ \,\text{and}\ \,G_{right}\le\Gamma\}.
$$

\section{Necessity in Theorem~\ref{040211a}}
Here we prove the necessity of Theorem~\ref{040211a}.
Throughout this section $\mZ_n$ is the additive group of integers modulo a positive integer~$n$.
\medskip

For any divisor $m$ of~$n$ denote by $i_{m,n}:\mZ_m\to\mZ_n$ and $\pi_{n,m}:\mZ_n\to\mZ_m$ the
group homomorphisms taking~$1$ to~$n/m$ and to~$1$ respectively. Using them we identify
the groups $i_{m,n}(\mZ_m)$ and $\mZ_n/\ker(\pi_{n,m})$
with~$\mZ_m$. Thus every section of $\mZ_n$ of order~$m$ is identified with the group~$\mZ_m$.
Moreover, the permutation $f\in\aut(\mZ_n)$ afforded by multiplication by an integer induces the
permutation $f^m\in\aut(\mZ_m)$ afforded  by multiplication by the same integer.\medskip

If~$\A$ is an S-ring over~$G=\mZ_n$ and $H$ is the $\A$-group of order~$m$, then $\A_H$ 
$\A_{G/H}$ are denoted respectively by $\A_m$ and~$\A^{n/m}$.
Let finally $\A_i$ be an S-ring over $\mZ_{n_i}$
($i=1,2$) and $(\A_1)^m=(\A_2)_m$ for some~$m$ dividing both~$n_1$ and $n_2$.
Then the unique S-ring~$\A$ over $\mZ_{n_1n_2/m}$ from~Theorem~A3.4
is denoted by $\A_1\wr_m\A_2$. We omit~$m$ if~$m=1$.\medskip 
%Given a group~$K\le\aut(G)$ and a
%prime~$p$ dividing $n$ we write $K_p$ for the $p$-projection of~$K$ in the sense
%of decomposition of the group $\aut(G)$ into the direct product of the groups $\aut(G_p)$ where
%$G_p$ is the Sylow $p$-subgroup of~$G$.

Below given a positive integer $m$ we set
$$
\Omega^*(m)=\css
\Omega(m),        &\text{if $m$ is odd,}\\
\Omega(m/2),      &\text{if $m$ is even}\\
\ecss
$$
where $\Omega(m)$ is the total number of prime factors of~$m$. We observe that $\Omega(m)\le 1$ if and
only if $m$ is a divisor of twice a prime number.

\thrml{040211b}
%Let $n=abcd$ where  such that $\GCD(ab,cd)=1$.
Let $n=n_1n_2$ where $n_1$ and $n_2$ are coprime positive integers such that
$\Omega^*(n_i)\ge 2$, $i=1,2$. Then a
cyclic group of order~$n$ is not Schur.
\ethrm
\proof Below for an integer $m\ge 3$ we denote by $K_m$ the subgroup of order~$2$
in the group $\aut(\mZ_m)$ that is generated by multiplication by $-1$. Suppose first that
$n_1=ab$ and $n_2=cd$ where $a,b,c,d\ge 3$ are integers. Set
\qtnl{240311z}
\A_1=\cyc(K_a\times K_c,\mZ_{ac}),\qquad
\A_2=\cyc(K_{bc},\mZ_{bc}),
\eqtn
\qtnl{240311y}
\A_3=\cyc(K_{ad},\mZ_{ad}),\qquad
\A_4=\cyc(K_{bd},\mZ_{bd}).
\eqtn
It is easily seen that the group $\aut(\A_i)$ is dihedral for $i=2,3,4$, and is the
direct product of two dihedral groups for $i=1$. Therefore the S-ring $\A_i$ is normal
for all~$i$. Moreover, $(\A_1)^c=\cyc(K_c,\mZ_c)=(\A_2)_c$ and $(\A_3)^d=\cyc(K_d,\mZ_d)=(\A_4)_d$.
Thus one can form S-rings
$$
\A_{1,2}=\A_1\wr_c\A_2\qquad\text{and}\qquad \A_{3,4}=\A_3\wr_d\A_4.
$$
It is easily seen that $(\A_{1,2})^{n_1}=\cyc(K_a,\mZ_a)\wr\cyc(K_b,\mZ_b)=(\A_{3,4})_{n_1}$.
Then
$$
\A:=\A_{1,2}\wr_{n_1}\A_{3,4}
$$
is an S-ring over $\mZ_n$. Thus it suffices to verify that $\A$ is not schurian.\medskip

Suppose on the contrary that $\A$ is schurian. Then 
%the schurity of the S-rings $\A_{1,2}$ and $\A_{3,4}$ and 
by Theorem~A1.2 %imply that 
the S-rings $\A_{1,2}$ and $\A_{3,4}$ are schurian, and
there exist groups
$\Delta_{1,2}\in\M(\A_{1,2})$ and $\Delta_{3,4}\in\M(\A_{3,4})$
such that
$$
(\Delta_{1,2})^S=(\Delta_{3,4})^S
$$
where $S$ is the section of order~$n_1$ used in the definition of the S-ring~$\A$.
In particular, for any permutation $f_1\in\Delta_{1,2}$ fixing~$0$ there exists a permutation 
$f_2\in\Delta_{3,4}$ fixing~$0$ and such that $f_1^S=f_2^S$. We claim: {\it the permutation $(f_1)^H$
where $H$ is the group of order~$ac$, is induced by multiplication by $\varepsilon\in\{1,-1\}$}.
However, if this is true, then the stabilizer of~$0$ in the group $(\Delta_{1,2})^H$
is contained in $K_{ac}$. Therefore the basic set of the S-ring associated with the former group
that contains $1$ is of cardinality $\le 2$. On the other hand, this S-ring coincides with $\A_1$
by the schurity of the S-ring $\A_{1,2}$ and the 2-equivalence of the groups $\Delta_{1,2}$ and
$\aut(\A_{1,2})$. So the above basic set has cardinality~$4$. Contradiction.\medskip

To prove the claim let $f_{1,1}$ and $f_{1,2}$ be the automorphisms of the S-rings $\A_1$ and $\A_2$ 
induced by~$f_1$, and $f_{2,3}$ and $f_{2,4}$ the automorphisms of the S-rings $\A_3$ and $\A_4$ 
induced by~$f_2$. Then the normality of these S-rings implies that $f_{1,1}\in K_a\times K_c$,
$f_{1,2}\in K_{bc}$, and that $f_{2,3}\in K_{ad}$ and $f_{2,4}\in K_{bd}$. Clearly,
\qtnl{220311z}
(f_{1,1})^c=(f_{1,2})^c\quad\text{and}\quad (f_{2,3})^d=(f_{2,4})^d
\eqtn
and due to the equality $(f_1)^S=(f_2)^S$ also
\qtnl{220311y}
(f_{1,1})^a=(f_{2,3})^a\quad\text{and}\quad (f_{1,2})^b=(f_{2,4})^b.
\eqtn
Next, the permutations $f_{1,2}$, $f_{2,3}$ and $f_{2,4}$ are induced respectively
by multiplications by some integers $\varepsilon_{1,2},\varepsilon_{2,3},\varepsilon_{2,4}\in\{1,-1\}$.
Therefore, by the second equalities of~\eqref{220311z} and~\eqref{220311y} we have
$$
\varepsilon_{1,2}=\varepsilon_{2,3}=\varepsilon_{2,4}.
$$
Denote this number by $\varepsilon$. Then by the first equalities of~\eqref{220311z}
and~\eqref{220311y} the permutations $(f_{1,1})^a$ and $(f_{1,1})^c$,
and hence the permutation $(f_1)^H$, are induced by multiplication by~$\varepsilon$.\medskip

To complete the proof we observe that the theorem is proved in all cases except for
the case when one of the numbers $n_1$, $n_2$, say $n_1$, is equal to~$8$. Then obviously
$n_1=ab/2$ and $n_2=cd$ where $a=b=4$ and $c,d\ge 3$ are odd integers. Let us define
S-rings $\A_1$, $\A_2$, $\A_3$ and $\A_4$ by formulas~\eqref{240311z}
and~\eqref{240311y}. Then again all these rings are normal, 
$$
(\A_1)^{2c}=\cyc(K_{2c},\mZ_{2c})=(\A_2)_{2c},\quad
(\A_3)^{2d}=\cyc(K_{2d},\mZ_{2d})=(\A_4)_{2d},
$$ 
and one can form S-rings
$$
\A_{1,2}=\A_1\wr_{2c}\A_2\qquad\text{and}\qquad \A_{3,4}=\A_3\wr_{2d}\A_4.
$$
It should be stressed that $\A_{1,2}$ and $\A_{3,4}$ are S-rings over the groups
$\mZ_{cn_1}$ and~$\mZ_{dn_1}$. 
It is also  easily seen that
$$
(\A_{1,2})^{n_1}=\cyc(K_a,\mZ_a)\wr_2\cyc(K_b,\mZ_b)=(\A_{3,4})_{n_1}.
$$
Then
$$
\A:=\A_{1,2}\wr_{n_1}\A_{3,4}
$$
is an S-ring over $\mZ_n$. Thus it suffices to verify that $\A$ is not schurian.
The rest of the proof repeats the proof of the first part literally.\bull\medskip

To complete the proof of the necessity we note that the required statement immediately follows from 
Theorem~\ref{040211b} and the lemma below.

\lmml{250311a}
An integer $n$ belongs to none of the families listed in~\eqref{250311b} if and only if $n=n_1n_2$
for some coprime positive integers $n_1$ and $n_2$ such that $\Omega^*(n_i)\ge 2$, $i=1,2$.
\elmm
\proof The sufficiency is straightforward by exaustive search. To prove the necessity let an integer
$n=p_1^{k_1}\cdots p_s^{k_s}$ belong to none of families~\eqref{250311b} where
$p_1,\ldots,p_s$ are pairwise distinct primes. Then without loss
of generality we can assume that
$$
2\le s\le 4\quad\text{and}\quad k_1\ge k_2\ge\cdots\ge k_s.
$$
Suppose on the contrary that $n$ cannot be decomposed into the product of coprime positive
integers $n_1$ and $n_2$ such that $\Omega^*(n_i)\ge 2$, $i=1,2$. Then $k_2=1$, for otherwise
$s=2$ and $n$ belongs to the third family with $p=2$, which is impossible. Thus
$k_2=\cdots=k_s=1$. Therefore $s=3$ or $s=4$, for otherwise $s=2$ and $n$ belongs to the second family.
Let $s=3$. Then $k_1\ne 1$ because $n$ does not belong to the fourth family.
So $k_1\ge 2$, and hence $2\in\{p_2,p_3\}$ by the supposition. However, then
$n$ belongs to the third family. Contradiction. Finally, let $s=4$. Then the supposition
implies that $k_1=1$ and one of the $p_i$'s equals~$2$. But then $n$ belongs to the fifth family.
Contradiction.\bull

\section{Quasidense S-rings}\label{130712a}
A circulant S-ring $\A$ is called {\it quasidense}, if any primitive
$\A$-section is of prime order. Any dense S-ring is obviously quasidense. It is
also clear that the property to be quasidense is preserved
by taking the restriction to any $\A$-section. Moreover, in the quasidense case
any minimal $\A$-group is of prime
order, any maximal $\A$-group is of prime index, and the S-ring $\A_S$ is
dense for any $\A$-section $S$ of prime power order.

\thrml{170112a}
Any quasidense circulant S-ring with trivial radical is cyclotomic, and hence dense.
\ethrm
\proof Let $\A$ be a quasidense circulant S-ring with trivial radical. Then $\A$ is the tensor product 
of a normal S-ring with trivial radical and S-rings of rank~$2$ by  
Theorem~A4.1. However, any normal circulant S-ring is cyclotomic by 
Theorem~A4.2. Besides, by the quasidensity the underlying group of any factor of rank~$2$ is 
of prime order. Therefore such a factor is also cyclotomic. Thus $\A$ is cyclotomic as the tensor 
product of cyclotomic S-rings.\bull\medskip

The following two statements will be used in proving the sufficiency of Theorem~\ref{040211a}
to find nontrivial $\A$-groups.

\crllrl{270112g}
Let $\A$ be a quasidense S-ring over a cyclic group~$G$. Then any subgroup of $G$ that
contains $\rad(\A)$ is an $\A$-group.
%$$\{H:\ \rad(\A)\le H\le G\}\subset\H(\A).$$
%any group $H$ such that $\rad(\A)\le H\le G$ is an $\A$-group. 
In particular, if $\rad(\A)_p=1$ for a prime divisor $p$ of $|G|$, then $G_{p'}$ is
an $\A$-group.  
\ecrllr
\proof The S-ring $\A_{G/L}$ where $L=\rad(\A)$, has trivial radical. Therefore by Theorem~\ref{170112a} it is 
dense. Thus required statement follows from the fact that a group $H$ containing $L$ is an $\A$-group if 
and only if the group $H/L$ is an $\A_{G/L}$-group.\bull

\crllrl{150612a}
Let $\A$ be a quasidense S-ring over a cyclic group~$G$. Suppose that $\A$ is not the $U/L$-wreath product
where $U/L$ is an $\A$-section such that the number $p:=|L|$ is prime. Then 
\nmrt
\tm{1} there exists $H\in\H(\A)$ such that $H\not\le U$ and $H_{p'}\in\H(\A)$,
\tm{2} if $q:=|G/U|$ is a prime other than $p$, then $H_{p'}\ge G_q$ for any group $H$ from statement~(1).
\enmrt
\ecrllr
\proof To prove statement~(1) we observe that by the hypothesis there exists 
$X\in\S(\A)$ outside $U$ such that $\rad(X)_p=1$. Then $H=\lg X\rg$ is an $\A$-group. 
Therefore the required statement follows from Corollary~\ref{270112g} applied
to the S-ring~$\A_H$. Next, the condition of statement~(2) implies that any group $H\not\le U$
contains a generator of $G_q$. Therefore $H\ge G_q$ which proves this statement.\bull\medskip

The following theorem reduces the schurity problem for circulant S-rings to the quasidense case. The proof
is based on the extension construction studied in~\cite{EP11}.

\thrml{220911a}
Given a circulant S-ring $\A$ there exists a quasidense S-ring $\A'\ge\A$ such
that $\A$ and $\A'$ are schurian or not simultaneously.
\ethrm
\proof Let us define an S-ring $\A'$
recursively as follows. If $\A$ has no singular class of composite order,
then we set $\A'=\A$; otherwise we set 
$$\A'=(\Ext_C(\A,\mZ S))'$$ where $C$ is
a singular class of composite order and $S=S_{min}(C)$. Then the S-ring $\A'$
has no singular classes of composite order. Moreover, from Theorem~A6.7
it follows that $\A$ and $\A'$ are schurian or not simultaneously. To complete the proof
let us verify that the S-ring $\A'$ is quasidense. Suppose on the contrary that this is not
true.  Then there exists a primitive
$\A'$-section $S$ of composite order. Then by Theorem~A4.6 the class
of projectively equivalent $\A'$-sections that contains $S$, is singular.
Contradiction.\bull\medskip

In general, the automorphism group of a quasidense S-ring is not solvable. However, from Theorem~A8.1
it follows that in the schurian case such an S-ring can always be obtained from an appropriate solvable 
permutation group in a standard way. The following theorem shows that ``locally''\ this group has a rather 
simple form. 

\thrml{230312a}
Let $\A$ be a schurian quasidense circulant S-ring. Then there exists a group $\Gamma\in\M(\A)$ such that
$\Gamma^S=\hol_\A(S)$ for any $\A$-section $S$ with $\rad(\A_S)=1$.
%\le\hol(S)$ and the group $\Gamma^S$ depends only on the S-ring~$\A_S$.
\ethrm

\rmrk
In fact, we prove that the equality in the theorem statement holds for any $S$ such that $\A_S$
is the tensor product of a normal S-ring and S-rings of rank~$2$. 
%In our case, by Lemma~\ref{010811a} the latter means that $\A_S$ is a cyclotomic S-ring with $|\rad(\A_S)|\le 2$.
\ermrk

\proof The quasidensity of~$\A$ implies that each primitive $\A$-section is of prime order. Therefore by 
Theorems~A4.6 and~A8.1 there exists a group $\Gamma\in\M(\A)$ such that $\Gamma^T\le\hol(T)$
for all primitive $\A$-sections~$T$. Let $S$ be an $\A$-section with $\rad(\A_S)=1$. Then
by Theorem~A4.1 the S-ring~$\A_S$ is the tensor product of a normal S-ring, say $\A_{T_0}$, 
and S-rings of rank~$2$, say $\A_{T_1},\ldots,\A_{T_k}$ where
$T_i$'s are $\A_S$-groups. It follows that 
$$
\Gamma^{T_0}\le\aut(\A_{T_0})\le\hol(T_0).
$$
Moreover, by the above $\Gamma^{T_i}\le\hol(T_i)$ for all~$i>0$, because the sections $T_1,\ldots,T_k$
are primitive. Thus
$$
\Gamma^S\le\prod_{i=0}^k \Gamma^{T_i}\le\prod_{i=0}^k\hol(T_i)=\hol(S).
$$ 
But then $\Gamma^S$ is obviously a unique subgroup of $\hol(S)$ in the set~$\M(\A_S)$.
Thus $\Gamma^S=\hol_\A(S)$.\bull

\section{Excluding families 1, 2 and 4}\label{190312a}

In the end of this section we prove the following theorem showing that any minimal non-schurian quasidense 
S-ring~$\A$ over a cyclic group $G$ contains two distinct minimal $\A$-groups and two distinct maximal 
$\A$-groups, the relationship between which, is as in Fig.~1.

\begin{figure}[h]
$\hspace{45mm}\xymatrix@R=10pt@C=20pt@M=0pt@L=5pt{
  &  & \VRT{G} \ar@{-}[dl]\ar@{-}[dr]                     &  & \\
  &  \VRT{U} \ar@{-}[dr] & & \VRT{V} \ar@{-}[dl]          & \\
  &  &  \VRT{} \ar@{--}[dd]                               &  & \\
  &  &                                                    &  & \\
  &  &  \VRT{} \ar@{-}[dl]\ar@{-}[dr]                     &  & \\
  &   \VRT{K} \ar@{-}[dr]  &  &  \VRT{L} \ar@{-}[dl]      & \\
  &  &  \VRT{1}                                               &  & \\
}$
\caption{}\label{f8}
\end{figure}

\thrml{270112b}
Let $\A$ be a minimal non-schurian quasidense S-ring over a cyclic group~$G$ of order~$n$ 
belonging to one of five families~\eqref{250311b}. Then 
\nmrt
\tm{1} $n$ belongs to the third or fifth families,
\tm{2} there exist distinct $\A$-groups $K$, $L$ of prime orders and distinct $\A$-groups $V$, $U$ of prime 
indices such that $LK\le U\cap V$ and $\A$ is a proper $U/L$-wreath product.
\enmrt
\ethrm

We begin with studying minimal non-schurian circulant S-rings. In the following statement we establish 
general properties of them.

\lmml{160112a}
Let $\A$ be a minimal non-schurian S-ring over a cyclic group~$G$. Then
\nmrt
\tm{1} $\A$ is a proper generalized wreath product,
\tm{2} if $\A$ is a proper $U/L$-wreath product, then 
%$\A_{U/L}$ is not normal and
%$\A_{U/L}$ is not the tensor product  of a normal S-ring and S-rings of rank~$2$; in particular, 
$\rad(\A_{U/L})\ne 1$; moreover, $|\rad(\A_{U/L})|>2$
whenever $\A_{U/L}$ is cyclotomic.
\enmrt
\elmm
\proof By Corollary~A4.3 we can assume  that $\rad(\A)\ne 1$. So statement~(1) follows
from Theorem~A4.1. The first part of statement~(2) follows from
the minimality of~$\A$, Theorem~A1.3 and Corollary~A1.4. Similarly, to prove the second part of this
statement it suffices to verify that the S-ring $\A_{U/L}$ is the tensor product of a normal 
S-ring and S-rings of rank~$2$ whenever it is cyclotomic and its radical is of order~$2$.
However, under
this condition the criterion of normality~\cite[Theorem~6.1]{EP01ce} implies that $\A_{U/L}$ is not normal 
if and only if it is the tensor product one factor of which is an S-ring of rank~$2$ (over a cyclic
group of prime order). Thus the required statement follows by induction.
\bull\medskip

By statement~(1) of Lemma~\ref{160112a} any minimal non-schurian circulant S-ring is a proper generalized 
wreath product. The following important theorem shows that in the quasidense case both operands are also 
proper generalized wreath products.

\thrml{160112b}
Let $\A$ be a minimal non-schurian quasidense S-ring over a cyclic group~$G$ of order~$n$
belonging to one of families~\eqref{250311b}. Then $\rad(\A_U)\ne 1$ and $\rad(\A_{G/L})\ne 1$
whenever $\A$ is a proper $U/L$-wreath product. 
\ethrm
\proof Let $\A$ be a proper $U/L$-wreath product. 
Suppose on the contrary that $\rad(\A_T)=1$ where $T\in\{U,G/L\}$.
Then the quasidensity of~$\A$ implies by Theorem~\ref{170112a} that the S-ring $\A_T$,
and hence the S-ring $\A_S$ with $S=U/L$,  is cyclotomic. By the minimality of~$\A$ and Lemma~\ref{160112a} 
this implies that $|\rad(A_S)|>2$. Thus $n$ does not belong to the fourth and the fifth families,
because otherwise $|S|$ is either prime, or $4$, or the product of 
two distinct primes. Moreover,
from Theorem~\ref{021109a} for $G=T$ it follows that
$S_l=1$ for some odd prime divisor $l$ of~$|T|$.
Thus $n$ does not belong to the first family. In the remaining two cases
the prime $l$ coincides with $p$, because otherwise $l=q$, and hence $|S|$ divides $2p$ which is impossible
by above. This proves the following lemma.

\lmml{070712a}
Under the above assumptions we have $n=pq^k$ or $n=2pq^k$, and $p\ne 2$. Moreover, 
\nmrt
\tm{1} if $\rad(\A_U)=1$, then $L_p\ne 1$,
\tm{2} if $\rad(\A_{G/L})=1$, then $(G/U)_p\ne 1$.\bull
\enmrt
\elmm

Let $\rad(\A_U)=1$. Assume that either $q=2$, or $G_{2'}$ is not an $\A$-group. 
By Lemma~\ref{070712a} the number $|G/L|$ is a power of~$q$ or twice a power of~$q$.
So by the assumption  there is a unique
maximal $\A_{G/L}$-group, say $U'/L$. Therefore $U'\ge U$, and hence $\A$ is the
$U'/L$-wreath product. Denote by $L'$ a maximal possible $\A$-group containing~$L$ for which
$\A$ is the $U'/L'$-wreath product. Then the uniqueness of $U'$ implies that 
$\rad(\A_{G/L'})=1$. (Indeed, otherwise by Theorem~A4.1 the S-ring 
$\A$ is the $U'/L''$-wreath product for some $L''>L'$, which contradicts the maximality of~$L'$).
Then by statement~(2) of Lemma~\ref{070712a} we conclude that $(G/U')_p\ne 1$.
Taking into account that $L'_p\ge L^{}_p\ne 1$, we conclude that $p^2$ divides $n$
which is impossible by Lemma~\ref{070712a}. This proves the first part of the following 
lemma (the second one is
proved in a similar way).

\lmml{050712s}
We have $n=2pq^k$ and $q\ne 2$. Moreover,
\nmrt
\tm{1} if $\rad(\A_U)=1$, then $G_{2'}$ is an $\A$-group,
\tm{2} if $\rad(\A_{G/L})=1$, then $G_2$ is an $\A$-group.\bull
\enmrt
\elmm

To complete the proof of Theorem~\ref{160112b} we come to a contradiction under the assumption
$T=U$ (the remaining case $T=G/L$ can be proved in a similar way). In this case we observe that
$U':=G_{2'}$ is an $\A$-group by Lemma~\ref{050712s}. We claim that
\qtnl{090712a}
\A=\A_{U'}\wr_{U'/L}\A_{G/L}.
\eqtn
Indeed, since $\A$ is the $U/L$-wreath product, it suffices to verify that $U'\ge U$. Suppose on the
contrary that this is not true. Then the number $|U|$ must be even. This implies that $G_2$ is an $\A$-group 
(we used the fact that the S-ring $\A_U$ is cyclotomic, and hence dense). Since $p\ne 2$ and $q\ne 2$, this
shows that $G_2$ is the $\A$-complement of~$U'$. Therefore by Corollary~\ref{030712a} we conclude that 
$\A=\A_{U'}\otimes\A_{G_2}$. By the minimality of~$\A$ this implies that the S-ring~$\A$ is schurian. The 
obtained contradiction proves equality~\eqref{090712a}.\medskip

%By  without loss of generality we can assume that the group $L$  is a 
%maximal $\A$-group such that $(G/L)_p=1$. 
After increasing the group $L$ in~\eqref{090712a} (if necessary), we can assume that it is a
maximal possible $\A$-group with that property. Then 
\qtnl{090712b}
(G/L)_p=1\qaq (G/L)_2\not\in\H(\A_{G/L}).
\eqtn
The first equality follows from Lemma~\ref{070712a}. To prove the second one suppose on the contrary
that the group $(G/L)_2$ is the $\A_{G/L}$-complement of~$U'/L$. Therefore by 
Corollary~\ref{030712a} we conclude that $\A_{G/L}=\A_{U'/L}\otimes\A_{(G/L)_2}$. By the minimality of~$\A$ 
and Theorem~\ref{140411a} this implies that the S-ring~$\A$ is schurian. The 
obtained contradiction proves~\eqref{090712b}.\medskip

Due to the quasidensity of~$\A$  formula~\eqref{090712b} implies that there is the only minimal 
$\A_{G/L}$-group, say $L'/L$, and $|L'/L|=q$. We claim that
\qtnl{090712c}
\A_{G/L}=\A_{U'/L}\wr_{U'/L'}\A_{G/L'}.
\eqtn
Indeed, otherwise by Corollary~\ref{150612a} 
%for $\A=\A_{G/L}$, $U=U'/L$, $L=L'/L$ and $p=q$
there exists an $\A_{G/L}$-group $H/L$ such that $H/L\not\le U'/L$ and $(H/L)_{q'}$ is an $\A_{G/L}$-group.
However, it is easily seen that in our case $(H/L)_{q'}=(G/L)_2$, which contradicts the second
relation in~\eqref{090712b}. The obtained contradiction proves~\eqref{090712c}.\medskip

Equalities \eqref{090712a} and~\eqref{090712c} show that the S-ring $\A$ is the $U'/L'$-wreath
product. However, this is impossible by the maximality of~$L$.\bull\medskip

{\bf Proof of Theorem~\ref{270112b}.} Statement~(1) immediately follows from statement~(2). To prove the 
latter we observe that by Lemma~\ref{160112a} the S-ring $\A$ is a proper $U/L$-wreath product
for some $\A$-groups $U$ and $L$. By the quasidensity of $\A$ we can assume that $L$ is of prime order
and $U$ is of prime index. Denote by
$\wt U$ a minimal subgoup of $U$ such that the S-ring $\A$ is the $\wt U/L$-wreath product. Then
by Theorem~\ref{160112b} the S-ring $\A_{\wt U}$ is a proper $U'/K$-wreath product
for some $\A$-groups $U'$ and $K$. Again we can assume that $K$ is of prime order. By the minimality 
of $\wt U$ we conclude that $K\ne L$. Besides, $KL\le U$ because
$L\le U$ and $K\le U'<\wt U\le U$.\medskip

 Similarly, denote by $\wt L$ a maximal subgroup of $U$ that contains~$L$ and such that the 
S-ring $\A$ is the $U/\wt L$-wreath product. Then again by Theorem~\ref{160112b} the S-ring $\A_{G/\wt L}$ is a 
proper $V/L'$-wreath product for some $\A$-groups $V$ and $L'$ such that 
$V$ is of prime index in~$G$. 
By the maximality of $\wt L$ we  conclude that $V\ne U$. Besides, obviously $V\ge L$.
To complete the proof it suffices to note that $K\le V$. Indeed, if this is not true, then
any nontrivial basic set of $\A_{G/\wt L}$ inside $K\wt L/\wt L$, is outside $V/\wt L$ 
and has trivial radical (because $|K\wt L/\wt L|=|K|$ is prime), which is impossible.\bull

\section{Excluding families 3 and 5 for $\mathbf{p\ne 2}$}\label{160312a}

In the end of this section we prove the following theorem that will enable us to exclude the
cases in the title.

\thrml{050612a}
Let $\A$ be a minimal non-schurian quasidense S-ring over a cyclic group~$G$ of order~$n$ 
belonging to the third or fifth of families~\eqref{250311b}. Then 
\nmrt
\tm{1} $p=2$,
\tm{2} $\A$ is both $U/L$- and $V/K$-wreath product where $K$, $L$, $U$, $V$ are $\A$-groups defined by
\nmrt
\tm{1}  $|K|=2$, $|L|=q$, $|U|=2qr$, $|V|=4q$ for $n=4qr$,
\tm{2}  $|K|=2$, $|L|=q$, $|U|=2q^k$, $|V|=4q^{k-1}$ for $n=4q^k$,
\enmrt
with $q$ and $r$ distinct odd primes and $k\ge 2$.
\enmrt
%Under the hypothesis of Theorem~\ref{270112a} suppose that $4$ does not divide $|U/L|$. Then
%$\A$ is the $V/K$-wreath product.
\ethrm

Throughout the rest of the section $\A$ denotes an S-ring satisfying the hypothesis of Theorem~\ref{050612a}.
It is also assumed that we are given $\A$-groups $K$, $L$, $U$, $V$ for which statement~(2) of
Theorem~\ref{270112b} holds.

\thrml{270112a}
The number $|U/L|$ is even.
\ethrm
\proof Suppose on the contrary that $|U/L|$ is odd. Then either $|L|=2$, or $|G/U|=2$. Let us consider
the former case, the latter one can be proved similarly. We are to find an $\A$-group $U'$ such that 
$L\le U'\le U$ and
\qtnl{270112d}
\A=\A_{U'}\wr_{U'/L}\A_{G/L}\quad\text{and}\quad (U')_{2'}\in\H(\A).
\eqtn 
Indeed, if such a group does exist, then by Corollary~\ref{030712a} with $G=U'$, $H=(U')_{2'}$, $S=L/1$ and 
$T=U'/L$, we obtain that $\A_U=\A_L\otimes\A_{U'/L}$. By the minimality of~$\A$ and Theorem~\ref{140411a} 
with $U=U'$ this implies that the S-ring $\A$ is schurian which is not true.\medskip

In the case $n=2pqr$ set $U'=U$. Then the left-hand side relation in~\eqref{270112d} is obvious. To
prove the other one we observe that by Theorem~A11.4 the number $|U|$ is the product of
three primes, the S-ring $\A_U$ is not a proper wreath product and $\A_{U/L}$ is a proper wreath
product. So, since the number $|U/L|$ is odd, the hypothesis of Lemma~A11.2 is satisfied 
for $\A=\A_U$, $S=U/L$ and $r=2$. By this lemma we obtain that $U_{2'}\in\H(\A_U)$, and we are done.\medskip

In the case $n=2pq^k$ set $U'$ to be a minimal $\A$-subgroup of~$U$, for which
the first relation in~\eqref{270112d} holds. We can assume that 
\qtnl{270112f}
(U')_p\ne 1.
\eqtn 
Indeed, otherwise $|U'|=2q^i$ for some $i$. Moreover, the minimality of $U'$ implies that the S-ring
$\A_{U'}$ is not the $U''/L$-wreath product where $U''$ is the subgroup of $U'$ of
index~$q$. Then by Corollary~\ref{150612a} with $p=2$ there exists an $\A_{U'}$-group 
$H_1$ such that $(U')_q\le H_1\le (U')_{2'}$. It follows that $H_1=(U')_q$, and 
the second relation in~\eqref{270112d} holds.\medskip

From~\eqref{270112f} it follows that $|G/V|=p$, and taking into account that $|U/L|$ is odd,
also that $p$ is odd. Therefore $U'\cap V\ne U'$. By the minimality of~$U'$ 
this implies that the S-ring $\A_{U'}$ is not a 
$(U'\cap V)/L$-wreath product. So by Corollary~\ref{150612a} with $(p,q)=(2,p)$ there exists an 
$\A$-group $H_1$ such that
\qtnl{270112z}
(U')_{p^{}}\le H_1\le (U')_{2'}.
\eqtn
Denote by $H$ a maximal
$\A$-subgroup of $U'$ that contains $H_1$. Then due to \eqref{270112z} we have $|U'/H|\in\{2,q\}$. If
$|U'/H|=2$, then  the second relation in~\eqref{270112d} holds and we are done. Finally, if
$|U'/H|=q$ then due to the minimality of~$U'$ the S-ring $\A_{U'}$,  is not an $H/L$-wreath product.
By Corollary~\ref{150612a} with $(p,q)=(2,q)$ there exists an $\A$-group $H$ such that
\qtnl{270112y}
(U')_{q^{}}\le H_2\le (U')_{2'}. %   ,\quad |H_2|\ \text{is odd},\quad H_2\in\H(\A).
\eqtn
Thus by~\eqref{270112z} 
and~\eqref{270112y} we have $(U')_{2'}=H_1H_2$, and hence $(U')_{2'}\in\H(\A)$.\bull

\thrml{170212a}
The order of $G$ is divisible by~$4$.
\ethrm
\proof Suppose first that $n=2pqr$. Then by statement~(1) of Theorem~A11.4 (where the
lattice of $\A$-groups is found) there are exactly two maximal $\A$-groups and exactly two minimal 
$\A$-groups; the former are of prime index whereas the latter are of prime order. From
statement~(2) of Theorem~\ref{270112b} it follows that these groups are $U$, $V$ and $K$, $L$ respectively,
and also that
\qtnl{170212b}
|U/L|\cdot|V/K|=|G|.
\eqtn
On the other hand, we claim that the S-rings $\A_U$ and $\A_{G/L}$ are non-normal. Indeed, otherwise by statement~(3)
of Theorem~A11.4 the number $|U/L|$ is a square of a prime. In our case this is possible only
for $p=2$. But in this case $|U/L|=4$ which is impossible by the latter theorem. The claim is proved.
So by statement~(5) of the same theorem the S-ring $\A$ is the $V/K$-wreath product. Then by Theorem~\ref{270112a}
of this paper the number $|V/K|$ is even. Thus $4$ divides $|G|$.\medskip

Let $n=2pq^k$. Suppose on the contrary that $p\ne 2$. Then since $q$ is odd and $|U/L|$ is even 
(Theorem~\ref{270112a}), from statement~(2) of Theorem~\ref{270112b} it follows that $|G/V|=2$ or
$|K|=2$. Let us consider the former case, the latter one can be proved similarly. In this case
$|V|$ is odd. Therefore by Theorem~\ref{270112a} the S-ring $\A$ is neither $V/K$- nor $V/L$-wreath product. 
So by Corollary~\ref{150612a} with $(p,q)=(2,p)$ and $(p,q)=(2,q)$ there exist $\A$-groups 
$H_1$ and $H_2$ such that 
$$
G_2\le H_1\le G_{p'}\qaq G_2\le H_2\le G_{q'}.
$$
Thus $G_2=H_1\cap H_2$, and hence $G_2$ is an $\A$-group. By Corollary~\ref{030712a} this
implies that $\A=\A_{G_2}\otimes\A_V$ which is impossible by the minimality of~$\A$.\bull

\thrml{210212a}
Without loss of generality we can assume that $4$ does not divide $|U/L|$.
\ethrm
\proof By Theorem~\ref{170212a} we have $p=2$. So $n=4qr$ or $n=4q^k$. In the former case 
$|U/L|$ is divisible by~$4$ only if $|U/L|=4$. However, in this case the S-ring $\A_{U/L}$
is cyclotomic and $|\rad(\A_{U/L})|\le 2$, which contradicts statement~(2) of Lemma~\ref{160112a}.
Thus we can assume that $n=4q^k$ and $4$ divides $|U/L|$. Then from statement~(2) of Theorem~\ref{270112b} it follows that 
$$
|G/U|=|L|=q\quad\text{and}\quad |G/V|=|K|=2.
$$ 
Therefore it suffices to verify that the S-ring $\A$ is either $U/K$- or $V/L$-wreath product. Suppose on
the contrary that this is not true. Then by Corollary~\ref{150612a} with $(p,q)=(2,q)$ and $(p,q)=(q,2)$ there 
exist $\A$-groups  $H_1$ and $H_2$ such that 
$$
G_q\le H_1\le G_{2'}\qaq G_2\le H_2\le G_{q'}.
$$
It follows that $H_1=G_q$ and $H_2=G_2$. Thus $G_q$ and $G_2$ are $\A$-groups. By the quasidensity of~$\A$ this 
implies that $\A$ is dense.\medskip

Denote by $\wt U$ the minimal $\A$-subgroup of $U$, for which the S-ring $\A$ is the $\wt U/L$-wreath product.
Then by Theorem~\ref{160112b} with $U=\wt U$ the radical of the ring $\A_{\wt U}$ is nontrivial.  Since this
S-ring is dense, from \cite[Theorem~3.4]{LM96} (see also statement~(1) of \cite[Theorem 5.4]{EP01ce})
it follows that it is a $U'/L'$-wreath product where the number 
$|L'|=|\wt U/U'|$ is the greatest prime divisor of $|\rad(\A_{\wt U})|$. By the minimality of 
the group $\wt U$ we conclude that this prime divisor is equal to~$2$. Thus
\qtnl{270211b}
L'=K=\rad(\A_{\wt U})\quad\text{and}\quad |\wt U/U'|=2.
\eqtn
By Corollary~\ref{030712a} we have $\A_{U'}=\A_K\otimes\A_{U'_q}$ and $\A_{\wt U/K}=\A_{\wt U_2/K}\otimes\A_{U'/K}$.
Since $\A_{\wt U}$ is the $U'/K$-wreath product this implies that 
\qtnl{270212c}
\A_{\wt U}=\A_{\wt U_2}\otimes\A_{\wt U_q}.
\eqtn
Therefore $\A_{\wt U/L}\cong\A_{\wt U_2}\otimes\A_{\wt U_q/L}$. Moreover, the S-ring $\A_{\wt U_2}$ being a
dense S-ring over a cyclic group of order~$4$, is cyclotomic and $|\rad(\A_{\wt U_2})|\le 2$.
Finally, by Corollary~\ref{210612a} the S-ring $\A_{\wt U_q/L}$ has trivial radical because by~\eqref{270211b} 
and~\eqref{270212c} so is the S-ring $\A_{\wt U_q}$. Therefore the latter S-ring is cyclotomic
by Theorem~\ref{170112a}. Thus the S-ring $\A_{\wt U/L}$ is cyclotomic and $|\rad(\A_{\wt U/L})|\le 2$
which contradicts statement~(2) of Lemma~\ref{160112a}.\bull\medskip

{\bf Proof of Theorem~\ref{050612a}.} By Theorem~\ref{270112b} we can assume that the hypothesis under
which Theorems~\ref{270112a}, \ref{170212a} and~\ref{210212a} were proved, holds for the S-ring~$\A$. 
Then statement~(1) immediately follows from Theorem~\ref{170212a}. Thus
$$
n=4qr\quad\text{or}\quad n=4q^k
$$
where $q$ and $r$ distinct odd primes and $k\ge 2$.  To prove statement~(2) choose the groups 
$K$, $L$, $U$, $V$ as above. Then obviously $|KL|=2q$ when $n=4q^k$. The same is also true for $n=4qr$
after interchanging $q$ and $r$ (if necessary). By Theorems~\ref{270112a} and~\ref{210212a}
we can also assume $|U/L|=2\,(\hspace{-2mm}\mod 4)$.\medskip 

Let $n=4qr$. Then by the above assumptions the number  $|U/L|$ is not a prime square. By statement~(3)
of Theorem~A11.4 this implies that neither of the S-rings $\A_U$ and $\A_{G/L}$ is normal. Therefore by 
statement~(5) of that theorem the S-ring $\A$ is the $V/K$-wreath product. Besides 
$|V/K|=2\,(\hspace{-2mm}\mod 4)$. Thus without loss of generality we can assume that $|K|=2$. Then $|L|=q$ and 
$|G/U|=2$. It follows that $|U|=2qr$ and $|V|=4q$, which completes the proof in this case.\medskip 

Let $n=4q^k$. Then by the above assumptions one of the following holds:
\nmrt
\tm{1} $|K|=2$, $|L|=q$, $|U|=2q^k$, $|V|=4q^{k-1}$,
\tm{2} $|K|=q$, $|L|=2$, $|U|=4q^{k-1}$, $|V|=2q^k$.
\enmrt
Thus it suffices to verify that $\A$ is the $V/K$-wreath product. Suppose that this is not true. 
%Then in each of these two cases we come to a contradiction separately.
\medskip

In case~(1) 
by Corollary~\ref{150612a} with $p=2$ there exists an $\A$-group~$H_1$ such that 
$G_q\le H_1\le G_{2'}$. It follows that $H_1=G_q$, and hence $G_q$ is an $\A$-group. By 
Corollary~\ref{030712a} this implies that 
$$
\A_U=\A_{G_q}\otimes\A_{K_{}}.
$$
Denote by $H$ the maximal $\A$-group such that $L\le H\le G_q$ and  $\rad(\A_H)=1$. Set $U'=HK$.
Then the radical of any basic set inside $U\setminus U'$ contains~$L$. Since the same is true 
also for any basic set outside~$U$, the S-ring $\A$ is the $U'/L$-wreath product. Besides,
$\rad(\A_{U'})=1$ because $\A_{U'}=\A_{H^{}}\otimes\A_{K^{}}$ and $|K|=2$. By statement~(2) of
Theorem~\ref{021109a} with $G=U'$ this implies that $\rad(\A_{U'/L})=1$. However,
this contradicts statement~(2) of Lemma~\ref{160112a}. Thus case~(1) is impossible.\medskip 

In case~(2) one can similarly prove that $G_2$ is an $\A$-group, and $\A_{G/L}=\A_{G_2/L}\otimes\A_{V/L}$.
It is also easily seen that any basic set outside $U$ is highest in $\A$ or in $\A_V$. Therefore
the S-ring $\A$ is the $U/L'$-wreath product where $L'=\rad(\A_V)$. Besides,
$\rad(\A_{G/L'})=1$ because $\A_{G/L'}=\A_{H/L'}\otimes\A_{V/L'}$ and $|H/L'|=2$,
where $H=G_2L'$. 
By statement~(2) of Theorem~\ref{021109a} with $G=G/L'$ this implies that $\rad(\A_{U/L'})=1$. 
However, this is impossible by statement~(2) of Lemma~\ref{160112a}.\bull

\section{Proof of Theorem~\ref{130112a}}\label{190312b}

Let $\A$ be a minimal non-schurian S-ring over a cyclic group~$G$ of order~$n$. Suppose on the contrary
that $n$ belongs to one of families~\eqref{250311b}. Since any divisor of $n$ also belongs to one of 
these families, by statement~(3) of Lemma~\ref{160112a} without loss of generality
we can assume that $\A$ is quasidense. Then by Theorems~\ref{270112b} and~\ref{050612a}
we have $n=4qr$ or $n=4q^k$, and $\A$ is both $U/L$- and $V/K$-wreath product where $K$, $L$, $U$, $V$ are 
$\A$-groups defined by
\nmrt
\tm{1}  $|K|=2$, $|L|=q$, $|U|=2qr$, $|V|=4q$ for $n=4qr$,
\tm{2}  $|K|=2$, $|L|=q$, $|U|=2q^k$, $|V|=4q^{k-1}$ for $n=4q^k$,
\enmrt
with $q$ and $r$ distinct odd primes and $k\ge 2$. In both cases we will verify that the hypothesis of 
Theorem~\ref{200312a} is  satisfied for some $\A$-groups so that the generalized wreath product for $\A$
defined there is proper. Then by that theorem $\A$ is schurian because due to the minimality of~$\A$ so are 
the operands of this product. Contradiction.\medskip

Suppose that we are in case~(1). Set $H_1=H_2=H$ where 
$H=KL=U\cap V$. First, we observe that relations~\eqref{080612a} and~\eqref{140512f} are obviously
satisfied. Furthermore, $|U/L|=2r$ is not a prime square, and hence by statement~(3) of Theorem~A11.4 the 
S-rings $\A_U$ and 
$\A_{G/L}$ are not normal. By statement~(4) of that theorem this implies that these S-rings 
are respectively the $H/K$- and $V/H$-wreath products. Besides, condition~(1) 
of Theorem~\ref{200312a} is trivially satisfied because the underlying groups of the S-rings 
$\A_{H/K}$, $\A_{U/H}$, $\A_{V/H}$ and $\A_{H/L}$,  are of prime orders, whereas condition~(2) 
is satisfied because $|K|=|G/U|=2$. Thus
the hypothesis of  Theorem~\ref{200312a} is satisfied.\medskip

Suppose that we are in case (2). To define the $\A$-groups from the hypothesis of Theorem~\ref{200312a}
we have to do preliminary work.  Set $M$ to be the minimal $\A$-subgroup of $G$ that contains $G_2$, 
and $N$ to be the maximal $\A$-subgroup of $G_q$. We claim that
\qtnl{170512g}
G_2\ne M,\quad  M_q\le N,\quad N\ne G_q.
\eqtn
Indeed, if $G_2=M$, then the radical of the highest basic set in $G_2$ has trivial $q$-part. However,
this is impossible because $\A$ is the $U/L$-wreath product. Similarly, if $N=G_q$, then the radical of 
the highest basic set in $G_q$ has trivial $2$-part. However,
this is impossible because $\A$ is the $V/K$-wreath product.
To prove the rest we observe that by Theorem~\ref{270112a} the S-ring~$\A$ is not a $U/K$-wreath product. 
Then by statement~(1) of Corollary~\ref{150612a} with $p=2$ there  exists an $\A$-group $H$ such that 
$G_2\le H$ and $H_q$ is an $\A$-group. So $M_q\le H_q$ and $N\ge H_q$
by the choice of~$M$ and $N$ respectively. This proves the claim.\medskip

Let us verify that the hypothesis of Theorem~\ref{200312a} is satisfied for $\A$-groups $K,\wt L,M,N,U,\wt V$
where
$$
\wt L=M\cap N\quad\text{and}\quad \wt V=MN.
$$
Then obviously $\wt L\le N$ and $M\le\wt V$. Moreover, from \eqref{170512g} it also follows that $H_1\le H_2$
where $H_1=K\wt L$ and $H_2=\wt V\cap U$. Since also 
\qtnl{040912a}
H_2=K\times N\quad\text{and}\quad G/H_1=M/H_1\times U/H_1,
\eqtn
the relations~\eqref{080612a} and~\eqref{140512f} hold. A part of 
the $\A$-group lattice is given at Fig.~\ref{f19}.\medskip

\def\VRT#1{*=<6mm>[o][F-]{#1}}
\begin{figure}[h]
$\hspace{35mm}
\xymatrix@R=10pt@C=20pt@M=0pt@L=5pt{
 & & & & \VRT{G} \ar@{-}[dl]\ar@{-}[dr]             &   \\
 & & &  \VRT{V} \ar@{--}[dl]\ar@{-}[dr] & & \VRT{U} \ar@{-}[dl]     \\
% & & \VRT{\wt V}  \ar@{--}[dl]\ar@{-}[dr]& & \VRT{}\ar@{--}[dl] & \\
 & & \VRT{\wt V}  \ar@{--}[dl]\ar@{-}[dr]& & \VRT{}\ar@{--}[dl] & \\
 & \VRT{M}\ar@{-}[dr] & & \VRT{H_2}\ar@{--}[dl]\ar@{-}[dr] & & \\
 & & \VRT{H_1} \ar@{-}[dr]\ar@{--}[dl]& & \VRT{N}\ar@{--}[dl] \\
 & \VRT{} \ar@{-}[dr]\ar@{-}[dl]& & \VRT{\wt L} \ar@{--}[dl]& & \\
 \VRT{K}\ar@{-}[dr] & & \VRT{L}\ar@{-}[dl] & & & \\
 & \VRT{1}   & & &  &   \\
}$
\caption{}\label{f19}
\end{figure}

To verify the rest of the hypothesis of Theorem~\ref{200312a} we observe that
condition~(2) is satisfied because $|K|=|G/U|=2$. We claim that
\qtnl{080612f}
\A=\A_U\wr_{U/\wt L}\A_{G/\wt L}.
\eqtn
Suppose on the contrary that this is not true. Then there exists a basic set~$X$ outside $U$ such that
$\rad(X)_q<\wt L$. Then obviously $M':=G_2\rad(X)$ is a proper subgroup of~$M$ that contains~$G_2$, which
is an $\A$-group
by Corollary~\ref{270112g}. However, this contradicts the minimality of~$M$. Next, let us verify that
\qtnl{080612g}
\A_U=\A_{H_2}\wr_{H_2/K}\A_{U/K}\quad\text{and}\quad\A_{G/\wt L}=\A_{\wt V/\wt L}\wr_{\wt V/H_1}\A_{G/H_1}.
\eqtn
To prove the first equality suppose on the contrary that the S-ring~$\A_U$ is not the $H_2/K$-wreath product. 
Then by statement~(1) of Corollary~\ref{150612a} with $p=2$ there  exists an $\A_U$-group $H\not\le H_2$ 
such that $H_q$ is an $\A$-group. However, this is impossible by the maximality of~$N$. The second equality can
be proved in a similar way. Thus by Remark~\ref{220612a} we only have to prove that
\qtnl{080612e}
\rad(\A_{H_2/K})=1,\quad \rad(\A_{U/H_1})=1,\quad \rad(\A_{H_2/H_1})=1.
\eqtn 

We observe that the third equality follows from the first one and Corollary~\ref{210612a} for $G=H_2/K$.
To prove the first equality in~\eqref{080612e} suppose on the contrary that $\rad(\A_{H_2/K})>1$. To get a
contradiction we use the idea from the proof of case~(1) in Theorem~\ref{050612a}. First, we
observe that by Corollary~\ref{030712a} we have 
$$
\A_{H_2}=\A_N\otimes\A_K.
$$
Set $U'=N'K$ where $N'$ is the maximal $\A$-group such that $L\le N'\le N$ and $\rad(\A_{N'})=1$.
Then $N'<N$ by the above supposition and the fact that $\A_N\cong\A_{H_2/K}$. Next, let $X$ be a
basic set outside $U'$. Then $L\le\rad(X)$ for $X\subset G\setminus U$ because $\A$ 
is the $U/L$-wreath product and for $X\subset H_2\setminus U'$ by the definition of $U'$.
The same is also true for $X\subset U\setminus H_2$. Indeed, otherwise  set $Q=\lg X\rg/\rad(X)$
and $S$ to be the image of the section $H_2/K$ in $Q$. Then $\rad(\A_S)=1$ by Theorem~\ref{021109a} 
applied to the S-ring $\A_Q$ and the section $S$. On the other hand, $\rad(\A_S)>1$ because $\rad(\A_{H_2/K})>1$ and $\rad(X)\le K$.
Contradiction. Thus the S-ring $\A$ is the $U'/L$-wreath product. Besides,
$\rad(\A_{U'})=1$ because $\rad(\A_{N'})=1$ and $A_{U'}=\A_{N'}\otimes\A_{K^{}}$.
By  statement~(2) of Theorem~\ref{021109a} for $G=U'$ this implies that $\rad(\A_{U'/L})=1$. However,
this contradicts statement~(2) of Lemma~\ref{160112a}. The second equality in~\eqref{080612e} 
is proved similarly following the proof of case~(2) in Theorem~\ref{050612a}.\bull

%\section{Restriction to a section in cyclotomic rings}\label{230911a}

\section{Auxiliary statements on S-rings}\label{130712b}

Given an S-ring $\A$ over a group $G$ we define an {\it $\A$-complement} of an $\A$-group~$H$ to be an
$\A$-group $H'$ such that $G=H\times H'$. When the group $G$ is cyclic, the group $H'$ is obviously uniquely
determined.

\thrml{240112c}
Let $\A$ be an S-ring over a cyclic group $G$. Suppose that an $\A$-group $H$ has an $\A$-complement and
$\A_S=\mZ S$ where $S$ is an $\A$-section projectively equivalent to $G/H$. Then 
given an $\A$-section $T$ projectively equivalent to $H$
the S-rings~$\A$ and $\A_S\otimes\A_T$ are Cayley isomorphic.
\ethrm
\proof Denote by $H'$ the $\A$-complement of $H$. Then obviously $H'/1$ and $G/H$ are respectively the 
smallest and greatest $\A$-sections in the class of projectively equivalent $\A$-sections that 
contains~$G/H$. This implies that the section $S$ is projectively equivalent to (in fact, a multiple of)
$H'/1$. By Theorem~A3.2 the S-rings $\A_{S^{}}$ and $\A_{H'}$ as well as
$\A_T$ and $\A_H$  are Cayley isomorphic. Thus without loss
of generality we can assume that $S=H'/1$ and $T=H/1$. Then $\A_{H'}=\mZ H'$, and hence 
%$h'X\in\S(\A)$ for all $h'\in H'$ and $X\in\A_H$. Therefore 
$$
\rk(\A)=|H'|\rk(\A_H)=\rk(\A_{H'})\rk(\A_{H^{}}).
$$
 Since also $\A\ge\A_{H^{}}\otimes\A_{H'}$ 
by Lemma~A2.1, we have $\A=\A_{H^{}}\otimes\A_{H'}$.\bull

\crllrl{030712a}
Theorem~\ref{240112c} remains true with the condition $\A_S=\mZ S$ replaced by $|S|=2$.\bull
\ecrllr

Some parts of the following statement appeared in a number of papers. Here we formulate it in a more or less
general form because it is used throughout the paper several times.

\thrml{021109a}
Let $\A$ be a cyclotomic S-ring with trivial radical over a cyclic group $G$. Suppose that $S$ is
an $\A$-section such that $S_p\ne 1$ for any odd prime divisor $p$ of $|G|$. Then
\nmrt
\tm{1} $|\rad(\A_S)|\le 2$,
\tm{2} $|\rad(\A_S)|=1$ unless $|S_2|=4$.
\enmrt
\ethrm
\proof
By \cite[Lemma~3.5]{LM96} given a set $X\in\S(\A)$ with $\rad(X)=1$ and a prime $p$ such that $p^2$ divides 
$m=|\lg X\rg |$, we have $\rad(X^p)=1$ unless $p=2$ and $m=8m'$ with $m'$ odd. This shows that 
$\rad(\A_U)=1$ where $U$ is the subgroup of $G$ of index~$p$, unless $p=2$ and $|G|=8m'$ with $m'$ odd.
Since $\A_{G/L}\cong\A_U$ where $L$ is the  subgroup of $G$ of order~$p$, we have $\A_{G/L}=1$
under the same conditions. Thus recursively applying these
results we reduces the lemma to the case 
$$
|S_2|\le 4,\quad |G_2|\le 8,\quad S_{2'}=G_{2'}.
$$
However, from \cite[Proposition~3.1]{LM96} with $m=|G|$ and $l=|G_{2'}|$ it follows that $\rad(\A_{G_{2'}})=1$. 
On the other hand, $\rad(\A_S)\le\rad(\A_{S_{2^{}}})\rad(\A_{S_{2'}})$ because
$\A_S\ge \A_{S_{2^{}}}\otimes\A_{S_{2'}}$ (see Lemma~A2.2). Thus 
$\rad(\A_S)\le\rad(\A_{S_{2^{}}})$, and we are done.\bull\medskip

From Theorems A4.1, A4.2 and \ref{021109a} we immediately obtain the following useful result.

\crllrl{210612a}
Let $\A$ be an S-ring with trivial radical over a cyclic $p$-group, $p$ odd. Then $\rad(\A_S)=1$ for
any $\A$-section~$S$.\bull
\ecrllr

The following statement gives a necessary and sufficient condition for the schurity of an $U/L$-wreath
product when the section $U/L$ is one of two sections forming an isolated pair of sections in the corresponding 
S-ring (see Definition~A6.1).

\thrml{140411a}
Let $\A=\A_U\wr_{U/L}\A_{G/L}$ be an S-ring over a cyclic group $G$. Suppose
that either $\A_U\cong\A_L\otimes \A_{U/L}$ or $\A_{G/L}\cong A_{U/L}\otimes \A_{G/U}$.
Then the S-ring $\A$ is schurian if and only if so are the S-rings $\A_U$ and $\A_{G/L}$.
\ethrm
\proof The necessity is obvious because given an $\A$-section $S$ the S-ring $\A_S$ is
schurian whenever so is $\A$. Let us prove the sufficiency under the assumption
$\A_U\cong\A_L\otimes \A_{U/L}$ (the rest can be proved analogously). Denote by
$f:U\to L\times (U/L)$ the corresponding Cayley isomorphism. Then $L^f=L$ and $H^{f}=U/L$
for a uniquely determined $\A$-group $H$. It follows that
$\A_U=\A_L\otimes\A_H$. Set
$$
\Delta_0=\aut(\A_{G/L})\quad\text{and}\quad
\Delta_1=\aut(\A_L)\otimes\Delta_H.
$$
where $\Delta_H$ is the full $f^{U/L}$-preimage of the group $(\Delta_0)^{U/L}$ in the group
$\aut(\A_H)$. Clearly,
$$
(G/L)_{right}\le\Delta_0,\quad U_{right}\le\Delta_1,\quad (\Delta_0)^{U/L}=(\Delta_1)^{U/L}.
$$
Moreover, by the schurity of the S-ring $\A_{G/L}$ the latter group is $2$-equivalent to
the group $\aut(\A_{U/L})$. It follows that the groups $\Delta_H$ and $\aut(\A_H)$, and
are $2$-equivalent. So by the schurity of the S-ring $\A_U$ the groups
$\Delta_1$ and $\aut(\A_U)=\aut(\A_L)\otimes\aut(\A_H)$ are also $2$-equivalent.
Thus by  Theorem~A1.2 the S-ring $\A$ is schurian and we are done.\bull

\section{A special generalized wreath product}\label{130712c}

In this section under special conditions we prove a necessary and sufficient condition for a $U/L$-wreath 
product to be schurian when the restriction of it to~$U/L$ is also a generalized wreath product.
We start with the description of elements of the canonical generalized wreath product introduced in Definition~A5.3
\footnote{The group that was denoted there by $\Delta_U$ is denoted here by $\Delta_1$.}.\medskip 

Let $G$ be an abelian group and $L\le U\le G$. Suppose we are given 
groups $\Delta_0\le\sym(G/L)$ and $\Delta_1\le\sym(U)$ such that $U/L$ is both $\Delta_0$- and 
$\Delta_1$-section and
$$
(G/L)_{right}\le\Delta_0,\quad U_{right}\le\Delta_1,\quad (\Delta_0)^{U/L}=(\Delta_1)^{U/L}.
$$
Then an element of the the canonical generalized wreath product 
$$
\Gamma=\Delta_1\wr_{U/L}\Delta_0
$$
can  explicitly be described as follows. Let us fix
bijections $h_X\in (G_{right})^{U,X}$ where $X\in G/U$. Suppose we are given a permutation $f_0\in\Delta_0$ 
and a family $\{f_X\in\Delta_1:\ X\in G/U\}$ of permutations such that
\qtnl{020412e}
(f_X)^{U/L}=(h_{X^{}})^{U/L}\,f_0^{X/L}\,((h_{X'})^{U/L})^{-1}
\eqtn
for all $X\in G/U$ where $X'$ is the $U$-coset for which $X'/L=(X/L)^{f_0}$. Then obviously there exists
a uniquely determined permutation $f\in\sym(G)$ for which 
$$
f^{G/L}=f_0\quad\text{and}\quad f^X=(h_X)^{-1}f_Xh_{X'}
$$
for all $X\in G/U$. We stress that this permutation depends on the choice of the permutations $h_X$. Denote 
it by $\{f_X\}\wr_{U/L} f_0$. Then the definition of the generalized wreath product of permutation groups 
implies immediately that
\qtnl{290312h}
\Gamma=\{\{f_X\}\wr_{U/L} f_0:\ f_0\in\Delta_0,\ f_X\in\Delta_1\ \text{for all}\ X\in G/U\}.
\eqtn

Let us turn to the main theorem of this section. Let $\A$ be a quasidense S-ring over a cyclic group~$G$. 
Suppose we are given $\A$-groups $K,L,M,N,U,V$ such that $L\le N$, $M\le V$, 
\qtnl{080612a}
H_1:=KL\le U\cap V:=H_2
\eqtn
and also
\qtnl{140512f}
H_2=K\times N\qaq G/H_1=M/H_1\times U/H_1.
\eqtn
The corresponding part of the $\A$-group lattice is represented in Fig.~\ref{f9}.

\begin{figure}[h]
$\hspace{45mm}\xymatrix@R=10pt@C=20pt@M=0pt@L=5pt{
  &  & \VRT{G} \ar@{-}[dl]\ar@{-}[dr]                     &  & \\
  &  \VRT{V} \ar@{-}[dr] & & \VRT{U} \ar@{-}[dl]          & \\
&    &  \VRT{H_2} \ar@{--}[dd]                            &  & \\
  &  \VRT{M}\ar@{--}[uu]\ar@{--}[dr] &     & \VRT{N}\ar@{--}[ul]\ar@{--}[dd] & \\
  &  &  \VRT{H_1} \ar@{-}[dl]\ar@{-}[dr]                  &  &   \\
  &   \VRT{K} \ar@{-}[dr]  &  &  \VRT{L} \ar@{-}[dl]      & \\
  &  &  \VRT{1}                                           &  & \\
}$
\caption{}\label{f9}
\end{figure}

\thrml{200312a}
In the above notation suppose that the S-rings $\A$, $\A_U$ and $\A_{G/L}$ are respectively the
$U/L$-, $H_2/K$- and $V/H_1$-wreath products such that
\nmrt
\tm{1} the S-rings $\A_{H_2/K}$, $\A_{U/H_1}$ and $\A_{V/H_1}$, $\A_{H_2/L}$ are of trivial radicals,
\tm{2} $\A_K=\mZ K$ and $\A_{G/U}=\mZ G/U$.
\enmrt
Then the S-ring $\A$ is schurian if and only if so are the S-rings $\A_U$ and $\A_{G/L}$.
\ethrm

\rmrkl{220612a}
Equalities~\eqref{140512f} together with condition~(2) imply by Theorem~\ref{240112c} that
$\A_{V/H_1}\cong\A_{G/U}\otimes\A_{H_2/H_1}$ and
$\A_{H_2/L}\cong\A_K\otimes\A_{H_2/H_1}$. Thus in our case the second part of condition~(1) is equivalent
to the equality $\rad(\A_{H_2/H_1})=1$.
\ermrk

\proof The necessity is obvious. To prove the sufficiency suppose that the S-rings $\A_U$ and $\A_{G/L}$ 
are schurian. Set
$$
\Gamma_1=\hol_\A(U/H_1)\quad\text{and}\quad\Gamma_2=\hol_\A(H_2/L).
$$
Then obviously $(\Gamma_1)^{H_2/H_1}=(\Gamma_2)^{H_2/H_1}$. So one can define the generalized wreath 
product $\Delta=\Gamma_2\wr_{H_2/H_1}\Gamma_1$. Thus by Theorem~A1.2 to complete
the proof it suffices to find 
groups $\Delta_1\in\M(\A_U)$ and $\Delta_0\in\M(\A_{G/L})$ such that
\qtnl{110512b}
(\Delta_1)^{U/L}=\Delta=(\Delta_0)^{U/L}.
%\Gamma_2\wr_{H_2/H_1}\Gamma_1=(\Delta_0)^{U/L}.
\eqtn
To do this we observe that by Theorem~\ref{230312a} and due to condition~(1) there exist
groups $\Gamma_3\in\M(\A_{U/K})$ and $\Gamma_6\in\M(\A_{V/L})$  such that 
\qtnl{110512a}
(\Gamma_3)^{H_2/K}=\hol_\A(H_2/K), \quad (\Gamma_3)^{U/H_1}=\hol_\A(U/H_1),
\eqtn
\qtnl{140512a}
(\Gamma_6)^{V/H_1}=\hol_\A(V/H_1),\quad (\Gamma_6)^{H_2/L}=\hol_\A(H_2/L).
\eqtn
Set
\qtnl{110512i}
\Gamma_4=\hol_\A(H_2),\quad \Gamma_5=\hol_\A(G/H_1).
\eqtn
Then clearly $(\Gamma_4)^{H_2/K}=(\Gamma_3)^{H_2/K}$ and $(\Gamma_6)^{V/H_1}=(\Gamma_5)^{V/H_1}$. Therefore
one can define generalized wreath products
$$
\Delta_1=\Gamma_4\wr_{H_2/K}\Gamma_3\qaq \Delta_0=\Gamma_6\wr_{V/H_1}\Gamma_5.
$$

First, let us prove that $\Delta_1\in\M(\A_U)$ and $\Delta_0\in\M(\A_{G/L})$. Indeed, since 
$U_{right}\le\Delta_1$ and $(G/L)_{right}\le\Delta_0$, it suffices to verify that 
$$
\Delta_1\twoe\aut(\A_U)\qaq \Delta_0\twoe\aut(\A_{G/L}).
$$
In its turn, to prove these relations it suffices to verify by Corollary~A5.7 applied to 
the S-ring $\A_U$ and the groups $\Gamma_4$, $\Gamma_3$, and the S-ring $\A_{G/L}$ and the groups 
$\Gamma_6$, $\Gamma_5$,   that
$$
\Gamma_4\twoe\aut(\A_{H_2}),\ \Gamma_3\twoe\aut(\A_{U/K}),
\ 
\Gamma_6\twoe\aut(\A_{V/L}),\ \Gamma_5\twoe\aut(\A_{G/H_1}).
$$
However, the statements on $\Gamma_3$ and $\Gamma_6$ hold by the definition of these groups.
Next, the hypothesis $\A_K=\mZ K$ implies by Theorem~\ref{240112c} that $\A_{H_2}$
is the tensor product of the cyclotomic rings $\A_K$ and $\A_N\cong\A_{H_2/K}$ (the latter S-ring
is cyclotomic by Theorem~\ref{170112a}). Therefore the S-ring $\A_{H_2}$ is cyclotomic, and hence
the groups $\aut(\A_{H_2})$ and $\Gamma_4$ are $2$-equivalent. Similarly, one can
prove that the group $\Gamma_5$ is $2$-equivalent to the group $\aut(\A_{G/H_1})$.\medskip

To prove \eqref{110512b} we note that due to \eqref{110512a}, \eqref{140512a}  and \eqref{110512i} we have
\qtnl{300312w}
(\Gamma_3)^{U/H_1}=\Gamma_1=(\Gamma_5)^{U/H_1}\quad\text{and}\quad 
(\Gamma_4)^{H_2/L}=\Gamma_2=(\Gamma_6)^{H_2/L}.
\eqtn
Therefore both $(\Delta_1)^{U/L}$ and $(\Delta_0)^{U/L}$ are contained in the group~$\Delta$.
%$\Gamma_2\wr_{H_2/H_1}\Gamma_1$. 
To prove the converse inclusion we observe that due to~\eqref{140512f}
there is an isomorphism
\qtnl{150512q}
G/V\to U/H_2,\quad X\mapsto X\cap U=:Y.
\eqtn
In what follows the factor sets $X$ and $Y$  modulo~$L$ are denoted by
$\ov X$ and $\ov Y$ respectively. For each $X\in G/V$ we fix a bijection $h_X\in(G_{right})^{V,\,X}$
that takes $H_2$ to $Y$, and set 
\qtnl{150512g}
h_Y=(h_X)^Y,\quad h_{\ov X}=(h_X)^{\ov X},\quad h_{\ov Y}=(h_X)^{\ov Y}.
\eqtn
Then due to equality~\eqref{290312h} any element $\ov f\in\Delta$ % \Gamma_2\wr_{H_2/H_1}\Gamma_1$
can be written in the form
$$
\ov f=\{f_{\ov Y}\}\wr_{H_2/H_1}\ov f_0
$$ 
for some permutation $\ov f_0\in\Gamma_1$ and a family of permutations $f_{\ov Y}\in\Gamma_2$ 
where $Y\in U/H_2$, such that 
\qtnl{020412h}
(f_{\ov Y})^{H_2/H_1}=(h_{\ov Y})^{H_2/H_1}\,(\ov f_0)^{Y/H_1}\, ((h_{\ov{Y'}})^{H_2/H_1})^{-1}
\eqtn
for all $Y\in U/H_2$ where $Y'$ is the $H_2$-coset in $U$ for which 
$Y'/H_1=(Y/H_1)^{\ov f_0}$. In what follows we find some elements of the groups
$\Delta_1$ and $\Delta_0$ the restrictions of which to $U/L$ coincide with $\ov f$.\medskip

To find the required permutation in $\Delta_0$ we observe that $\A_{M/H_1}=\mZ M/H_1$
because $\A_{M/H_1}\cong\A_{G/U}$, and the latter is a group ring by condition~(2). So by Theorem~\ref{240112c} we have
\qtnl{150512w}
\Gamma_5=(M/H_1)_{right}\otimes\Gamma_1.
\eqtn 
Therefore this group contains the permutation 
$f_0=\id_{M/H_1}\otimes\ov f_0$. Clearly,
\qtnl{150512a}
(f_0)^{U/H_1}=\ov f_0.
\eqtn
Next, let $X\in G/V$. Then by~\eqref{140512a} there exists a permutation $f_{\ov X}\in\Gamma_6$ 
that leaves the set $H_2/L$ fixed and such that 
\qtnl{150512r}
(f_{\ov X})^{H_2/L}=f_{\ov Y}.
\eqtn 
Below we show that
\qtnl{150512b}
(f_{\ov X})^{V/H_1}=(h_{\ov X})^{V/H_1}\,(f_0)^{X/H_1}\, ((h_{\ov{X'}})^{V/H_1})^{-1}
\eqtn
where $X'$ is the $V$-coset in $G$ for which $X'/H_1=(X/H_1)^{f_0}$. Then
one can define a permutation $f=\{f_{\ov X}\}\wr_{V/H_1} f_0$ belonging to the group~$\Delta_0$. 
This is we wanted to find because $f^{U/L}=\ov f$ by~\eqref{150512g},
\eqref{150512a} and~\eqref{150512r}.\medskip

To prove~\eqref{150512b} we observe that 
$(\Gamma_6)^{V/H_1}=(V/H_2)_{right}\otimes\hol_\A(V/M)$
because $\A_{V/H_2}=\mZ V/H_2$ (see above). So
$$
(f_{\ov X})^{V/H_1}=(f_{\ov X})^{V/H_2}\otimes(f_{\ov X})^{V/M}=
(f_{\ov X})^{V/H_2}\otimes(f_{\ov Y})^{V/M}.
$$
On the other hand, the permutation $f_{\ov X}$ leaves the set
$H_2/L$ fixed. So $(f_{\ov X})^{V/H_2}$ lives the set $H_2$ fixed. Since also 
$(f_{\ov X})^{V/H_2}\in (V/H_2)_{right}$, this implies that $(f_{\ov X})^{V/H_2}=\id_{V/H_2}$.
Thus~\eqref{150512b} holds by \eqref{020412h} and the choice of the bijections $h_X$.\medskip

To find a permutation $f\in\Delta_1$ such that $f^{U/L}$ coincides with the permutation $\ov f$ 
defined in~\eqref{150512q}, we observe that due to~\eqref{110512a} there exists
a permutation $f_0\in\Gamma_3$ such that equality~\eqref{150512a} holds.
Next, for each $Y\in  U/H_2$ we define a permutation of $H_2/K$ defined by
\qtnl{020412a1}
g_Y=(h_Y)^{H_2/K}\,(f_0)^{Y/K}\,((h_{Y'})^{H_2/K})^{-1}
\eqtn
where $Y'$ is the $H_2$-coset in $U$ for which $Y'/K=(Y/K)^{f_0}$. However, the bijection $h_X$
by its choice leaves the set $U$ fixed. So 
$$
(h_Y)^{H_2/K}\,(f_0)^{Y/K}\,((h_{Y'})^{H_2/K})^{-1}=((h_X)^{U/K}\,f_0\,((h_{X'})^{U/K})^{-1})^{H_2/K}
$$
Thus $g_Y$ belongs to the group $(\Gamma_3)^{H_2/K}=(\Gamma_4)^{H_2/K}$. Therefore,
due to condition~(2) of the theorem the permutation
$$
f_Y:=g_Y\otimes (f_{\ov Y})^{H_2/N}
$$
belongs to the group $\Gamma_4$. Moreover, equalities \eqref{150512a}, \eqref{020412a1} and~\eqref{020412h}
imply that
$$
(f_Y)^{H_2/L}=(h_Y)^{H_2/L}\,(f_0)^{Y/L}\, (h_{Y'})^{H_2/L})^{-1}.
$$
Thus one can define a permutation $f=\{f_Y\}\wr_{H_2/L} f_0$ belonging to the group~$\Delta_1$. By
the choice of $f_0$ and $f_Y$ we have $f^{U/L}=\ov f$, which completes the proof.\bull

\end{document}